\documentclass[12pt,a4paper]{article}
\usepackage{amsmath}
\usepackage{amssymb}
\usepackage{graphicx}
\usepackage[left=1in,right=1in,top=1in,bottom=1in]{geometry}
\usepackage{cite}
\usepackage[british]{babel}

\newtheorem{theo}{Theorem}[section]
\newtheorem{lm}{Lemma}[section]

\newtheorem{rmk}{Remark}[section]

\newtheorem{proposition}{Proposition}[section]

\allowdisplaybreaks

\numberwithin{equation}{section}

\def\R{{\mathbb R}}
\def\Z{{\mathbb Z}}
\def\N{{\mathbb N}}
\def\S{{\mathcal S}}
\def\Exp{{\mathbf E}}

\def\Pr{{\mathbf P}}
\def\1{{\mathbf 1}}

\def\eps{\varepsilon}
\def\omu{{\overline \mu}}
\def\umu{{\underline \mu}}
\newcommand{\tod}{{\stackrel{d}{\longrightarrow}}}

\def\bx{{\mathbf x}}
\def\bu{{\mathbf u}}
\def\c{{\mathrm{c}}}

\newcommand{\essinf}{\mathop{{\rm ess~inf}}}
\newcommand{\esssup}{\mathop{{\rm ess~sup}}}

\newcommand{\proof}{\noindent {\bf Proof.}~}

\newcommand{\qed}{\hfill $\square$}

\newcommand{\F}{{\mathcal F}}

\title{Rate of escape and central limit theorem for the supercritical Lamperti problem}

\author{Mikhail V. Menshikov\footnote{Department of Mathematical Sciences, Durham University, South Road, Durham DH1 3LE, UK.} 
\and Andrew R. Wade\footnote{Department of Mathematics and Statistics, University of Strathclyde, 26 Richmond Street, Glasgow G1 1XH, UK.}}

\date{30 July 2010}

\begin{document}

\maketitle

\begin{abstract}
The study of discrete-time stochastic processes on the half-line
with mean drift at $x$ given by $\mu_1 (x) \to 0$ as $x \to \infty$
is known as {\em Lamperti's problem}.
We give sharp almost-sure bounds for
processes of this type in the case where $\mu_1 (x)$ is of order
$x^{-\beta}$ for some $\beta \in (0,1)$.
The bounds are of order $t^{1/(1+\beta)}$, so the process
is super-diffusive but sub-ballistic (has zero speed).
We make minimal assumptions on the moments of the increments
of the process (finiteness of $(2+2\beta+\eps)$-moments for our main results,
 so $4$th moments certainly suffice)
and do not assume that the process is time-homogeneous or
Markovian. In the case where $x^\beta \mu_1 (x)$
has a finite positive
limit, our results imply a strong law of large numbers, which strengthens and
generalizes earlier results of Lamperti and Voit. We prove
an accompanying central limit theorem, which appears
to be new even in the case of a nearest-neighbour random walk, although our
result is considerably more general. This answers
a question of Lamperti. We also prove transience of the process
under weaker conditions than those that we have previously seen
in the literature. Most of our results also cover the case where $\beta =0$.
We illustrate our results
with applications to birth-and-death chains and to multi-dimensional
non-homogeneous random walks.
 \end{abstract}



\section{Introduction}
\label{sec:intro}

 In a pioneering series of
 papers \cite{lamp1,lamp2,lamp3}
 published in the early 1960s,
 J.~Lamperti
 systematically studied
 how the asymptotic behaviour
 of a nonnegative real-valued discrete-time
 stochastic process
 with asymptotically zero drift is
 governed by the (first two)
 moment functions of its
 increments. In the last two decades there has been renewed interest
  in Lamperti's problem
 and in particular in its applications
 to studying the behaviour of
 complicated multidimensional
 processes (see e.g.~\cite{FMM,mvw}).
 A
 special case of Lamperti's problem
 supported on $\Z^+ := \{0,1,2,\ldots\}$
 is that of asymptotically-zero-drift
 birth-and-death chains,
 for which exact calculations
 are often possible (using for instance
 Karlin--McGregor theory \cite{csvd,km1,km2});
 although classically well-studied,
 there has been recent renewed interest
 in such birth-and-death chains (see e.g.~\cite{cfr1,cfr2}),
 particularly
 in the context of modelling random
 polymers (see e.g.~\cite{dec,giacomin}).
 The study of continuous-time analogues of the general Lamperti
 problem seems to have begun only recently: see e.g.~\cite{ds}.

 Let us describe informally Lamperti's problem.
 Consider a stochastic process $X = (X_t)_{t \in \Z^+}$
 on $[0,\infty)$. For now suppose that $X$ is a time-homogeneous Markov
 process (that is, a Markov process with stationary transition
 probabilities)
 and that its increment moment functions
 \begin{equation}
 \label{mu1}
  \mu_k (x) = \Exp [ (X_{t+1} - X_t)^k \mid X_t = x] \end{equation}
 are well-defined for $k \geq 0$;
 one way to ensure this is to impose a uniform bound
 on the increments. (We will relax all of these conditions shortly.)
  Lamperti's problem is
 to determine how the asymptotic behaviour of $X$
 depends upon $\mu_1$ and $\mu_2$.

Under mild regularity conditions,
 the behaviour of $X$
 is rather standard when, outside some bounded set,
  $\mu_1(x) \equiv 0$
 (the zero drift case)
or  $\mu_1(x)$ is  
  uniformly bounded
 to one side of $0$. Roughly speaking,
 in the zero-drift case $X$ behaves
 like a simple symmetric random walk
 and is null-recurrent, in the case
 of uniformly negative drift $X$
 is positive-recurrent with exponentially
 decaying
 stationary distribution,
 and in the case of uniformly
 positive drift $X$ is transient
 with positive speed (i.e., ballistic).

 This motivates the study of the {\em asymptotically zero drift}
 regime, in which $\mu_1(x) \to 0$
 as $x \to \infty$, 
  to investigate phase transitions.
  It turns out that there is a rich spectrum
  of possible behaviours of $X$, governed
  by $\mu_1$ and $\mu_2$; 
  we  mention heavy-tailed positive-recurrence,
  transience with sub-linear rate of escape
  (diffusive and super-diffusive motion both being possible),
  weak convergence to a  Bessel process, and so on.

Results of Lamperti \cite{lamp1,lamp3} imply that from the
point of view of the recurrence classification of $X$, the
case where $|\mu_1 (x)|$ is of order $x^{-1}$ and $\mu_2 (x)$
is of order 1 is critical. In the present paper we are interested
in the {\em supercritical} case where $\mu_1 (x)$ is positive and of order
$x^{-\beta}$, $\beta \in (0,1)$. Here, under mild conditions,
transience is assured: our primary interest is 
to quantify this transience by studying the {\em rate
of escape} and accompanying second-order behaviour.

  As well as being of interest in their own right, stochastic processes
on the half-line with mean
drift asymptotically zero are important for the study of multidimensional processes
by the method of Lyapunov-type functions (see e.g.~\cite{FMM}). In this
context it is particularly
desirable to work in some generality without
imposing, for instance, assumptions of the Markov
property, a countable state-space, or uniformly bounded increments.
Thus we work in  
more generality
 than the model outlined informally above.
 To start with, the assumption on uniformly
 bounded increments can   be relaxed, and
 replaced by an appropriate moments condition. Another important relaxation
 (building on the ideas in Lamperti's first paper on the topic \cite{lamp1})
 is that we do not need $X$ to be a Markov process. It
 is invaluable with regard to applications to be able to dispense with
the Markov assumption. The prototypical
 illustration of this
latter point is provided by the case where $X$ is given by $X_t = \| Y_t \|$,
the norm of some multidimensional (perhaps Markov but not necessarily spatially homogeneous)
process. If $Y_t$ has mean drift zero, $X_t$ will typically have
$\mu_1 (x) \to 0$ as $x \to \infty$.

Relaxing the Markov assumption leads to a slight
 complication in defining the correct analogues
of (\ref{mu1}), but does not complicate our proofs which
are based on general martingale arguments. The process $X$
will be taken to be adapted to some filtration $(\F_t)_{t \in \Z^+}$.
Important families of processes that fit into our framework
include non-Markov processes where $\F_t = \sigma (X_0, X_1, \ldots, X_t)$
and the law of $X_{t+1}$ depends on the entire previous history of the process,
as well as processes where $X_t$ is not Markov by itself, but $X_t =
f (Y_t)$ for some Markov process $Y_t$ on a general space $\Sigma$,
a measurable function
$f: \Sigma \to [0,\infty)$, and $\F_t = \sigma (Y_0,Y_1, \ldots, Y_t)$. The first of
these two situations was treated by Lamperti in \cite[Section 3]{lamp1}, and the second in
\cite[Section 4]{lamp1} (see also \cite[Section 5]{lamp3}); we work somewhat more generally.

 In the next section we will
 describe more precisely the model
 that we consider  and give our main results.
 In Section \ref{exam} we give two applications of our results
 to stochastic processes of interest in their own right.
 The first is the birth-and-death chain case; even in
 this classical setting, some of our
 results seem to be new. Our second example is
 a model inaccessible to many classical methods: a multi-dimensional
 non-homogeneous random walk. In the latter setting, our results
 add to the analysis of MacPhee {\em et al.} \cite{mmw}.

 \section{Model, results, and discussion}
 \label{sec:model}

\subsection{The model and main results}

 We now introduce our notation and assumptions.
  Let $X = (X_t)_{t \in \Z^+}$  be a discrete-time
 stochastic process adapted to a filtration
$(\F_t)_{t \in \Z^+}$
 and
taking values in an unbounded
subset $\S$ of $[0,\infty)$. In applications
$\S$ may be countable (e.g.~the birth-and-death chain example
 in Section \ref{exam1}) or uncountable
(e.g.~the non-homogeneous random walk
example in Section \ref{exam2}, or
the application to stochastic
billiards in \cite{mvw});
it is thus desirable to make no further restriction
on $\S$.
 
The central object in all that follows will
be the conditional mean increment (the one-step
mean drift)
$\Exp [ X_{t+1} - X_t \mid \F_t ]$.
Many of the conditions in our theorems will suppose that
an inequality hold involving the $\F_t$-measurable
random variables $\Exp [ X_{t+1} - X_t \mid \F_t ]$ and $X_t$;
such inequalities will have to hold
a.s.~and in an appropriate asymptotic sense (as   $X_t \to \infty$).
It will be convenient therefore to introduce some notation for
 upper and lower bounds
on the mean increment $\Exp [ X_{t+1} - X_t \mid \F_t ]$
as   functions of $X_t$. 

Shortly we will define 
$\umu_1:  \S \to \R$ and $\omu_1:   \S \to \R$  such that for all $t \in \Z^+$,  
\begin{equation}
\label{muprop}
 \umu_1 (X_t)   \leq \Exp [  X_{t+1} - X_t  \mid \F_t ]
\leq \omu_1 (X_t)  , ~{\rm a.s.}. \end{equation}
  If $X$ is a Markov process, $\Exp [  X_{t+1} - X_t  \mid \F_t ] =
 \Exp [  X_{t+1} - X_t  \mid X_t ]$, a.s., and we can   take
 \[ \umu_1 (x ) = \inf_{t \in \Z^+} \Exp [  X_{t+1} - X_t  \mid X_t =x ], ~~~{\rm and}~~~\omu_1 (x ) = \sup_{t \in \Z^+} \Exp [  X_{t+1} - X_t  \mid X_t =x ];\]
 if additionally $X$ is time-homogeneous then $\umu_1 (x) \equiv \omu_1 (x) \equiv \mu_1 (x)$ where
 $\mu_k : \S \to \R$ is given by
   \begin{equation}
    \label{mu2}
     \mu_k (x)
    = \Exp [ (X_{t+1} -X_t)^k \mid X_t = x ] ~~~(t \in \Z^+),\end{equation}
    provided the expectation exists.   Loosely speaking, in the general
    case $\Exp [  X_{t+1} - X_t  \mid \F_t ]$ involves  additional randomness in $\F_t$, once $X_t$
    has been fixed. Thus $\omu_1 (x)$ should be the (essential) supremum over
    this additional randomness given $\{ X_t = x\}$.   For $\umu_1$ the
    situation is analogous. 
    
    Let us now formally define $\umu_1$ and $\omu_1$.
    Suppose that $\Exp [  X_{t+1} - X_t  \mid \F_t ]$ exists for all $t \in \Z^+$.
    By standard theory of conditional expectations (see e.g.~\cite[Section 9.1]{chung}),
    for each $t \in \Z^+$
     there
    exist a Borel-measurable function $\phi_t : \S \to \R$
    and an $\F_t$-measurable random variable $\psi_t$ such that
    $\Exp [ \psi_t \mid X_t ] = 0$   and, a.s.,
    \begin{equation}
    \label{orthog}
     \Exp [  X_{t+1} - X_t  \mid \F_t ] = \Exp [ X_{t+1} - X_t \mid X_t ] + \psi_t = \phi_t (X_t ) + \psi_t .\end{equation}
Set $\mu_1 (t ; x) := \phi_t (x) + \psi_t$,
an $\F_t$-measurable random variable.
Then for $x \in \S$  
   define
\begin{align}
\label{mudef1}
 \omu_1 (x) & := \sup_{t \in \Z^+}
 \esssup \mu_1 (t ; x) , \\
 \label{mudef2}
 \umu_1 (x) & := \inf_{t \in \Z^+} 
 \essinf \mu_1 (t ; x) .  \end{align}
Provided the expectations in
question exist,
$\omu_1(x), \umu_1(x)$ are (non-random) $\R$-valued
functions of $x \in \S$; clearly
$\omu_1(x) \geq \umu_1(x)$ for all $x \in \S$.
Then (\ref{mudef1}) and (\ref{mudef2}) define
 functions with the property (\ref{muprop}). We 
 provide some further discussion of the
 definitions in (\ref{mudef1}) and (\ref{mudef2}),
 and   give some illustrative examples,
 in Section \ref{tech} below. 

In the time-homogeneous Markov case, the statement
  of our results is simplified,  and
  $\omu_1$, $\umu_1$ can be replaced simply by $\mu_1$
 defined by (\ref{mu2}) everywhere.  One such example, which might also
 be useful for orientation purposes, is the birth-and-death chain example
 described in Section \ref{exam1} below. 
  As mentioned above, in applications it
 can be important to dispense with the Markovian
 assumption. It often turns out to be the case
 in applications that as $x \to \infty$, $\omu_1 (x) \sim \umu_1 (x)$;
 the  example in Section \ref{exam2} below  demonstrates such a case,
 and also the importance of not having
 to assume a Markov property for $X$.

 Returning to the general setting,
  for our purposes the most interesting case is when
 $\omu_1 (x), \umu_1 (x) \to 0$ as $x \to \infty$.
 Results of Lamperti \cite{lamp1,lamp3} show that from the
 point of view of the recurrence classification
 of $X$, the case where $\omu_1 (x), \umu_1 (x)$ are of order
 $1/x$ is critical (assuming some natural
 regularity conditions). Our focus in the present paper
 is the {\em supercritical} case where
 $\omu_1 (x), \umu_1 (x)$ are of order $x^{-\beta}$
 (in the positive direction)
 for some $\beta \in (0,1)$. In this case Lamperti
 \cite{lamp1} proved that $X$ is transient
 (that is, $X_t \to \infty$ a.s.)
 under certain regularity assumptions;
 we give a proof of this result under weaker
 conditions (Theorem \ref{trans}). Our primary interest, however,
 is the nature of the transience, in particular
 the {\em rate of escape}, i.e., the speed at which $X_t \to \infty$.
 The results of this paper give sharp bounds of order $t^{1/(1+\beta)}$
 for $X_t$ (Theorem \ref{thm1}), which in the special case where $\omu_1 (x) \sim \umu_1 (x) \sim \rho x^{-\beta}$
 imply a strong law of large numbers (Theorem \ref{llnthm})
 that improves upon results
 of Lamperti \cite{lamp2} and Voit \cite{voit92}. We also study
 the second-order behaviour, obtaining a central limit theorem (Theorem \ref{cltthm})
 to accompany the law of large numbers. Although not our primary concern,
most of our results also cover the case where $\beta =0$.

Let us  state
our basic assumption:
\begin{itemize}
\item[(A0)] Let $X=(X_t)_{t \in \Z^+}$ be
a stochastic process on the unbounded set $\S \subseteq [0,\infty)$
adapted to the filtration $(\F_t)_{t\in\Z^+}$.
Suppose that for some $x_0 \in \S$, $\Pr[ X_0 \leq x_0 \mid \F_0 ] =1$.
\end{itemize}
We also assume the following condition:
\begin{itemize}
\item[(A1)]
Suppose that
$\limsup_{t \to\infty} X_t = \infty$ a.s..
\end{itemize}
Condition (A1) 
is necessary for our questions of interest to be non-trivial,
and is usually straightforward to verify in a particular application: for instance, a sufficient condition is that for any $y \in (0,\infty)$
there
exist $w:\Z^+  \to \Z^+$ and $\eps > 0$ such that
\[ \inf_{t \in \Z^+} \Pr [ X_{t+w(t)} > y \mid \F_t ] > \eps , {\rm a.s.}.\]
Indeed, if $X$ is an irreducible time-homogeneous
 Markov chain
and $\S$ is countable, (A1) holds
automatically. For suitable concepts of irreducibility
in more general state-spaces, see \cite{mt}.

We also need to assume some regularity condition on the increments
of $X$. For our purposes, we will need a moment
bound of the form
\begin{equation}
\label{moments}
  \sup_{t \in \Z^+}  \Exp [  | X_{t+1} -X_t |^{\gamma} \mid \F_t ] \leq B , ~{\rm a.s.},  \end{equation}
for some $B < \infty$ and $\gamma > 0$.
If (\ref{moments}) holds with $\gamma \geq 1$,
 $\omu_1$ and $\umu_1$ 
 given by (\ref{mudef1}) and (\ref{mudef2})
exist as
$\R$-valued functions.
Assumption of  (\ref{moments})  amounts to, in some sense, the choice
of a correct scale for the process $X$.

Our first result yields transience of the supercritical Lamperti problem.

\begin{theo}
\label{trans}
Suppose that (A0) and (A1) hold, and that there exists $\beta \in [0,1)$  such that
(\ref{moments}) holds for some $\gamma > 1+\beta$ and
\[ \liminf_{x \to \infty}  ( x^\beta \umu_1 (x) ) >0  .\]
Then $X$ is transient, i.e.,
$X_t \to \infty$ a.s.~as $t \to \infty$.
\end{theo}

Theorem \ref{trans} proves transience  under weaker conditions than we have seen
 previously published;
for instance Lamperti \cite[Theorem 3.2]{lamp1}
(see also \cite[Section 9.5.3]{mt})
assumed (\ref{moments}) with $\gamma >2$ and also
that $\Exp [ (X_{t+1} -X_t)^2 \mid \F_t ] \geq v$ a.s.~for $v >0$;
Lamperti \cite{lamp1} was mainly concerned
with the critical case ($\beta=1$), where such stronger conditions are natural,
but they are not necessary here, as Theorem \ref{trans} shows.

Next we move on to our main topic, the quantitative asymptotic
behaviour of $X$.
The first natural question is what bounds we can
obtain under conditions of comparable strength
to those in Theorem \ref{trans}. We have the
following upper bound.

\begin{theo}
\label{thm0}
 Suppose that (A0) holds, there exists $\beta \in [0,1)$
 such that
 \[ \limsup_{x \to \infty}  ( x^\beta \omu_1(x) ) < \infty ,\]
 and   (\ref{moments}) holds for some $\gamma > 1+\beta$.
  Then for any $\eps>0$, a.s., for all but finitely many $t$,
 \begin{equation}
 \label{crudebound}
  \sup_{0 \leq s \leq t} X_s \leq t^{\frac{1}{1+\beta}} (\log t)^{\frac{1}{1+\beta}+\eps}.\end{equation}
\end{theo}

Next we impose stronger conditions on $X$ in order
to obtain a tighter upper bound, as well as a
complementary lower bound. Our bounds
will involve the constants $\lambda(a,\beta)$
defined
for $a \in (0,\infty)$,
$\beta \in (0,1)$ by
  \begin{equation}
  \label{ddef}
  \lambda(a,\beta ) := (  a (1+\beta) )^{\frac{1}{1+\beta}} .\end{equation}
 The next result gives sharp almost-sure bounds on $X$.

 \begin{theo}
 \label{thm1}
Suppose that (A0) and (A1) hold, and
 that for some $\beta \in [0,1)$ and some $a, A
 \in (0,\infty)$ with $a \leq A$,
 \begin{equation}
 \label{mubounds}
  a = \liminf_{x \to \infty} ( x^\beta \umu_1 (x) ) \leq
 \limsup_{x \to \infty} ( x^\beta \omu_1 (x) ) = A .\end{equation}
 Suppose that (\ref{moments}) holds for some $\gamma > 2+2\beta$.
 Then, a.s.,
 \[ \lambda(a,\beta) \leq
 \liminf_{t \to \infty}  \frac{X_t}{t^{1/(1+\beta)}} \leq \limsup_{t \to \infty}  \frac{X_t}{t^{1/(1+\beta)}}
 \leq \lambda(A,\beta) .\]
  \end{theo}

\begin{rmk}
The proof of the upper bound on $X_t$ given by Theorem \ref{thm1}
  only uses the condition on $\omu_1$ in (\ref{mubounds}) and not the
  condition on $\umu_1$ there.\end{rmk}

Note that since $\beta <1$, certainly taking $\gamma =4$ in (\ref{moments})
suffices for Theorem \ref{thm1}. 
Theorem \ref{thm1} implies that in the case  $\beta \in (0,1)$
the transience given in Theorem \ref{trans}
is {\em super-diffusive} but {\em sub-ballistic}, since $1/2 < 1/(1+\beta) < 1$.
This should be contrasted with the critically transient case ($\beta=1$)
where the drift is $O(x^{-1})$ and $X$ is transient,
in which
case there are upper and lower bounds for $X_t$ of order about $t^{1/2}$ known under
additional conditions, see \cite[Section 4.1]{mvw},
where for instance it is shown in \cite[Theorem 4.2]{mvw}
that $X_t \geq t^{1/2} (\log t)^{-D}$ for some $D \in (0,\infty)$
and all but finitely many $t$ (in the critically transient
birth-and-death chain case, certain
sharp  bounds are a byproduct of the invariance principle of \cite{cfr2}).

  An immediate corollary of Theorem \ref{thm1}, obtained on taking $a=A=\rho$, is the following
  strong law of large numbers.

 \begin{theo}
 \label{llnthm}
Suppose that (A0) and (A1) hold, and that for some $\beta \in [0,1)$,
 \begin{equation}
 \label{limits}
  \lim_{x \to \infty} x^{\beta} \omu_1 (x) = \lim_{x \to \infty} x^{\beta} \umu_1 (x) = \rho \in (0,\infty) .\end{equation}
 Suppose that (\ref{moments}) holds for some $\gamma > 2+2\beta$.
 Then as $t \to \infty$, a.s.,
 \begin{equation}
 \label{lln}
    \frac{X_t}{t^{1/(1+\beta)}} \longrightarrow \lambda(\rho, \beta) .\end{equation}
 \end{theo}

Lamperti   \cite[Theorem 7.1]{lamp2}
  obtained a weaker version of Theorem \ref{llnthm} under
  more restrictive conditions. 
  Specifically, \cite[Theorem 7.1]{lamp2} assumes
that  $X$ is a time-homogeneous Markov process with $\lim_{x \to \infty} x^{\beta} \mu_1 (x)  = \rho$
  and $\sup_x |\mu_k (x)| < \infty$ for all $k$,
  where $\mu_k$ is given by (\ref{mu2}). Then  \cite[Theorem 7.1]{lamp2} says that
  (\ref{lln}) holds
  {\em with convergence in probability}. Lamperti \cite[p.~768]{lamp2}
asks whether his result ``can be strengthened to almost sure convergence'';
Theorem \ref{llnthm} answers this affirmatively, and also shows that the assumptions
in \cite{lamp2} can be relaxed to a significant extent. Theorem \ref{llnthm} also generalizes
a result of Voit \cite{voit92} in the birth-and-death chain case: see Section \ref{exam1} below.

It is natural to ask whether,
under the assumptions of Theorem \ref{llnthm}, there is a central limit theorem
  to accompany the law of large numbers.
This question was raised by Lamperti \cite[p.~768]{lamp2}, and seems
to have remained open even for the case of a birth-and-death chain.
The following result shows that there is a central limit theorem,
provided that we impose a somewhat stronger version of (\ref{limits})
and an asymptotic stability condition on the second moments of the
increments. Here and subsequently
`$\tod$' denotes convergence in distribution. Unlike
our preceding results, the case $\beta =0$ is excluded from the
following theorem.

 \begin{theo}
 \label{cltthm}
Suppose that (A0) and (A1) hold, and that for some $\beta \in (0,1)$ and $\rho \in (0,\infty)$,
as $x \to \infty$,
 \begin{equation}
 \label{limits2}
\umu_1(x) = \rho x^{-\beta} + o ( x^{-\beta - \frac{1-\beta}{2}} ); ~~~
\omu_1(x) = \rho x^{-\beta} + o ( x^{-\beta - \frac{1-\beta}{2}} ).
\end{equation}
  Suppose that (\ref{moments}) holds for some $\gamma > 2+2\beta$, and
 that for some $\sigma^2 \in (0,\infty)$,
 \begin{equation}
 \label{var}
 \Exp [ (X_{t+1} -X_t )^2 \mid \F_t] \to \sigma^2 , ~{\rm a.s.},
 ~{\rm as}~t \to \infty.
 \end{equation}
 Then as $t \to \infty$, 
 \begin{equation}
 \nonumber
    \frac{X_t - \lambda (\rho, \beta) t^{1/(1+\beta)}}{t^{1/2}} \tod
    Z \sigma \sqrt{ \frac{1+\beta}{1+3\beta}} ,\end{equation}
    where $Z$ is a standard normal random variable.
 \end{theo}

\subsection{Further remarks on $\umu_1$ and $\omu_1$; examples}
\label{tech}
 
 We now briefly discuss further the definitions in (\ref{mudef1})
and (\ref{mudef2}), and give some examples
for particular classes of process $X$ that should help to clarify
the nature of the crucial functions $\umu_1$ and $\omu_1$.
Recall that
\[ \esssup \mu_1 ( t ; x) = \inf \{ z \in \R : \Pr [ \mu_1 ( t ; x ) > z ] = 0 \} ,\]
with a similar expression for $\essinf$.
 Some intuitive feeling for the quantities $\umu_1$, $\omu_1$ is best gained by
specializing our general framework to some particular families of processes. 

\paragraph{Markov processes} If $\F_t = \sigma (X_0,\ldots,X_t)$ and $X$ is Markov,
we have that $\Exp [ X_{t+1} - X_t \mid \F_t ] = \Exp [ X_{t+1} - X_t \mid X_t ]$, a.s., so that,
with the notation at (\ref{orthog}),
\[ \mu_1 (t ; x) = \phi_t (x) = \Exp [ X_{t+1} - X_t \mid X_t = x ], ~{\rm a.s.}, \]
for any $t$.
When the state-space $\S$ is countable, this last quantity is simply expressed in
terms of the one-step transition probabilities $\Pr [ X_{t+1} = y \mid X_t = x ]$.
In the case of general $\S$, $\mu_1 (t ;x)$ can be expressed in terms
of a corresponding Markov transition kernel. In either case, we then have
that $\omu_1 (x) = \sup_{t \in \Z^+} \Exp [ X_{t+1} - X_t \mid X_t = x]$,
with a similar expression for $\umu_1(x)$. If $X$ is additionally time-homogeneous,
$\Exp [ X_{t+1} - X_t \mid X_t = x ]$ does not depend on $t$ so that
$\omu_1 (x) \equiv \umu_1 (x)$.

\paragraph{History-dependent processes} Suppose, more generally, that
 $\F_t = \sigma (X_0,\ldots,X_t)$
and the law of $X_{t+1}$ depends only upon $(X_0,\ldots,X_t)$.
For convenience, take $\S$ to be countable. Then we can write
\[ \Exp [ X_{t+1} - X_t \mid \F_t ] = \sum_{x_0,\ldots, x_t \in \S} \Exp [
X_{t+1} - X_t \mid X_0 = x_0, \ldots,  X_t = x_t ] \1 \{ X_0 = x_0, \ldots,  X_t = x_t \}.\]
This last expression can be written as $\mu_1 (t ; X_t)$ where $\mu_1 (t ; x)$ is given by
\[    \sum_{x_0, \ldots, x_{t-1} \in \S}  \Exp [
X_{t+1} - X_t \mid X_0 = x_0, \ldots, X_{t-1} = x_{t-1}, X_t = x  ] \1 \{ X_0 = x_0, \ldots, X_{t-1} = x_{t-1} \}.\]
It follows that in this case
\[ \omu_1 (x) = \sup_{t \in \Z^+} \sup_{
\stackrel{x_0, \ldots, x_{t-1} \in \S:}{\Pr [ X_0 = x_0, \ldots, X_{t-1} = x_{t-1} ] > 0}}
 \Exp [
X_{t+1} - X_t \mid X_0 = x_0, \ldots, X_{t-1} = x_{t-1}, X_t = x  ] ,\]
with an analogous expression for $\umu_1$.
In the case where $\S$ is uncountable, the expressions are similar
but may be understood in terms of regular conditional distributions.
 This
 formulation is essentially used by Lamperti \cite[p.~322]{lamp1}.  

\paragraph{Functions of Markov processes} Suppose that $(Y_t)_{t \in \Z^+}$
is a Markov process on some state-space $\Sigma$, and that for some
measurable function
$f : \Sigma \to [0,\infty)$, $X_t = f(Y_t)$. Set
$\F_t = \sigma (Y_0, \ldots, Y_t)$. Then $X$ has state-space $\S = f (\Sigma)$, and 
$X$ is typically
non-Markovian; see e.g.~\cite{rosenblatt} for a discussion on the latter point. 
Now $\Exp [ X_{t+1} - X_t \mid \F_t ] =  \Exp [ X_{t+1} - X_t \mid Y_t ]$, a.s.,
and if $\Sigma$
 (hence $\S$) is countable, we may write
 \[ \Exp [ X_{t+1} - X_t \mid \F_t ] = \sum_{x \in \S} \sum_{y \in \Sigma : f(y) = x} \Exp [ X_{t+1} - X_t \mid Y_t = y ] \1 \{ Y_t =y , X_t = x \} .\]
Expressing the latter quantity as $\mu_1 (t; X_t)$ entails
 \[ \mu_1 (t ; x) = \sum_{y \in \Sigma : f(y) = x} \Exp [ X_{t+1} - X_t \mid Y_t = y ] \1 \{ Y_t =y  \}.\]
 It follows that, in this case,
 \[ \omu_1 (x) = \sup_{t \in \Z^+} \sup_{y \in \Sigma : f(y) =x ,\, \Pr [ Y_t =y ] >0 } \Exp [ X_{t+1} - X_t \mid Y_t = y ] ,\]
 and similarly for $\umu_1$.
This situation often arises in applications, where $f$ may be,
for instance, a
Lyapunov-type function applied to a multi-dimensional process.
See the example in Section \ref{exam2} below, as well
as \cite[Section 4]{lamp1} and \cite[Section 5]{lamp3}.

\subsection{Open problems and paper outline}

We finish this section by mentioning some possible directions for future work.
A natural question is whether Theorem \ref{thm1} holds
under a weaker moments condition.  Also of
interest is whether any weak limit theory
analogous to Theorem \ref{cltthm} is available
when (\ref{limits}) holds but (\ref{limits2}) does not.

In \cite{mv}, an analogue of Lamperti's problem was considered
for processes with $\Exp [ X_{t+1} - X_t \mid X_t = x] \approx c x^\alpha t^{-\beta}$,
loosely speaking. It seems likely that for appropriate $\alpha, \beta$ one could
obtain results similar to ours in that setting.

 The outline of the remainder of the paper is as follows.
 In Section \ref{exam} we discuss two applications of our main theorems,
 specifically to birth-and-death chains (nearest-neighbour
 random walks on $\Z^+$) in Section \ref{exam1},
 and to non-homogeneous random walks in $\R^d$ in Section \ref{exam2}.
 Section \ref{sec:proofs} is devoted to the proofs of
 our theorems. In Section \ref{outline} we give a brief
 overview of our proofs.
 In Sections \ref{proof0}, \ref{proof1}, \ref{cltproof} and \ref{proof2}
  we prove Theorems \ref{thm0}, \ref{thm1}, 
  \ref{cltthm}
  and  \ref{trans} respectively; finally in Section \ref{rwproof}
 we prove our result (Theorem \ref{rwthm}) on non-homogeneous random walk
 presented in Section \ref{exam2} below.

 \section{Applications}
 \label{exam}

 \subsection{Birth-and-death chains}
 \label{exam1}

 Suppose that $X$ is an irreducible time-homogeneous Markov chain
 supported on the countable set $\S = \Z^+$
with jumps of size at most $1$.
  Specifically, suppose that
there exist  sequences $a_x, b_x, c_x$ ($x \in \N := \{1,2,3,\ldots\}$)
with $a_x >0$, $b_x \geq 0$, $c_x >0$ and $a_x + b_x + c_x =1$ for all $x \in \N$.
 Define the transition law of $X$ for $t \in \Z^+$ as follows: for $x \in \N$,
\begin{align*}
\Pr [ X_{t+1} = x+1 \mid X_t = x ]& = a_x , \\
\Pr [ X_{t+1} = x \mid X_t = x ] & = b_x, \\
\Pr [ X_{t+1} = x-1 \mid X_t = x ] & = c_x,
\end{align*}
and with reflection from $0$ governed by
$\Pr [ X_{t+1} = 1 \mid X_t =0 ] = 1$. Of course
in this setting $X$ has uniformly bounded
increments, so that (\ref{moments}) holds for all $\gamma>0$,
 and is an irreducible time-homogeneous
Markov chain on $\Z^+$, so that (A1) holds as well.

Such an $X$ is known as a {\em birth-and-death
chain} or birth-and-death random walk.
Such processes have been extensively studied in various contexts, and
are often
amenable to explicit computation.
Early contributions to the theory of such random walks, 
particularly to the recurrence/transience
classification, are due to Harris
\cite{harris} and Hodges and Rosenblatt \cite{hr}. Orthogonal
polynomials provide 
one fruitful tool for analysis of such processes (see e.g.~\cite{csvd} for a survey); 
this approach 
dates back at least to Karlin and McGregor \cite{km1,km2}.

For $x \in \N$, with $\mu_1(x)$ and $\mu_2(x)$ defined
by (\ref{mu2}) we have that
\[ \mu_1(x) =  a_x - c_x , ~~~
\mu_2(x) = 1 - b_x > 0.\]
The asymptotically zero drift case is the case where
$
\lim_{x \to \infty} (a_x - c_x) =0$.
There is an extensive literature concerned
with various special cases where
$| x \mu_1 (x) | = O(1)$. For recent
work, we refer to \cite{dec,cfr1,cfr2};
the papers of Cs\'aki, F\"oldes and R\'ev\'esz
cited include references to some of the older literature.
 We are in the supercitical
case if, for $\beta \in (0,1)$,
\begin{equation}
\label{sup1}
\lim_{x \to \infty} x^\beta (a_x - c_x ) = \rho \in (0,\infty).\end{equation}
In this case the following law of large numbers
 is due to Voit \cite[Theorem 2.11]{voit92}
(in fact Voit works in a more general setting of random walks on polynomial
hypergroups, which do not concern us here).
Note that there is a misprint
in the limiting constant in the 
statement of Theorem 2.11 of \cite{voit92} (the proof there does yield the
 correct constant): the $1/(1+\alpha)$ power should
be applied to the entire limiting expression, not just the $\mu$ there;
this typo persists into \cite[Theorem D]{cfr1}.

 \begin{proposition}
 \label{voitprop}
 \cite[Theorem 2.11]{voit92}
 Suppose that $X$ is a birth-and-death chain specified by $a_x, b_x, c_x$ as
 described above. Suppose that (\ref{sup1}) holds for $\beta \in (0,1)$
 and $\rho \in (0,\infty)$. Suppose also that the two limits
 \begin{equation}
 \label{voit}
 \lim_{x \to \infty} a_x  ~{\textrm{and}} ~ \lim_{x \to \infty} c_x ~\textrm{exist in}~ (0,1), \end{equation}
 (in which case they must take the same value). Then
(\ref{lln}) holds.
\end{proposition}

  Proposition \ref{voitprop} is a special case of
   our Theorem \ref{llnthm}, under the additional
   assumption (\ref{voit}). Theorem \ref{llnthm} shows that
   the assumption (\ref{voit}) is not necessary for the result: only the mean
   is important, not the absolute probabilities of going left or right.
   
Our Theorem \ref{cltthm} above has the following immediate
(and apparently new) consequence in the birth-and-death chain case. 
Under the assumption that $\lim_{x \to \infty} b_x =b$, (\ref{voit}) holds (with
limit $\frac{1-b}{2}$ for $a_x$ and $c_x$),
so this central limit theorem can be seen as the natural
companion to Voit's law of large numbers \cite[Theorem 2.11]{voit92}.

 \begin{theo}
 \label{clt2}
 Suppose that $X$ is a birth-and-death chain specified by $a_x, b_x, c_x$ as
 described above. Suppose that for some $\beta \in (0,1)$ and $\rho \in (0,\infty)$,
 \[ a_x - c_x = \rho x^{-\beta} + o ( x^{-\beta - \frac{1-\beta}{2}} ); ~~~ \lim_{x \to \infty} b_x = b \in [0,1) .\] 
  Then as $t \to \infty$, 
  for a standard normal random variable $Z$,
 \begin{equation}
 \label{clt3}
    \frac{X_t - \lambda (\rho, \beta) t^{1/(1+\beta)}}{t^{1/2}} \tod
    Z  \sqrt{ \frac{(1-b)(1+\beta)}{1+3\beta}} .\end{equation}
 \end{theo}

We make some final remarks on the case, of secondary interest to us here, where $\beta=0$.
   Our Theorem \ref{llnthm} applies to the case $\beta=0$, i.e., where
   $a_x - c_x \to \rho \in (0,1]$ as $x \to \infty$, in which case our result says that
   $t^{-1} X_t \to \rho$ a.s.~as $t \to \infty$. This particular result
   has been previously obtained by Pakes \cite[Proposition 4]{pakes},
   under some   more restrictive conditions, including $\rho = 1$ and $b_x \equiv 0$,
   and also, for general $\rho$ but again under conditions more restrictive than ours,
   in a result of Voit \cite[Corollary 2.6]{voit90}. 
   In the case $\beta =0$, the second-order behaviour of $X$ is somewhat different
  (our Theorem \ref{clt2} does not apply). See for instance \cite[Theorem 7]{pakes}
  and \cite[Theorems 2.7--2.10]{voit90}.

 \subsection{Rate of escape for non-homogeneous random walk on $\R^d$}
  \label{exam2}
  
  In this section we illustrate the application of our 
  results to a non-homogeneous random walk model
  similar to that of \cite{mmw}. Fix
  $d \in \{2,3,\ldots\}$.
  Let $\Xi = (\xi_t)_{t \in \Z^+}$
  be a time-homogeneous Markov process with state-space an unbounded
  subset $\Sigma$ of $\R^d$. The law of the increment $\xi_{t+1} - \xi_t$ then
  depends only on the position of $\xi_t$; this is formalized
  in general in terms of Markov transition kernels
  (cf.~\cite[Section 3.4]{mt}), so that we may use the
  notation $\Pr [ \, \cdot \mid \xi_t = \bx]$ for the conditional
  distributions and $\Exp [ \, \cdot \mid \xi_t = \bx]$ for the
  corresponding expectations.
   
Write $\| \, \cdot \, \|$ for
  the Euclidean norm on $\R^d$ and ${\bf 0}$ for the origin.
    We use the notation $\bx$ for a point of $\R^d$,
  and, when $\bx \neq {\bf 0}$,
   $\hat \bx := \bx / \| \bx \|$ for the corresponding unit
  vector. We use `$\cdot$' to denote the usual scalar product
  on $\R^d$.
   We assume that
  there exist $\rho \in (0,\infty)$ and $\beta \in (0,1)$ 
   such that, for $\bx \in \Sigma$,
  \begin{equation}
  \label{ass1}
  \Exp [ (\xi_{t+1} - \xi_t ) \cdot \hat \bx  \mid \xi_t = \bx ] = \rho \| \bx \|^{-\beta} + o (\| \bx \|^{-\beta} ) ,
  \end{equation}
  as $\| \bx\| \to \infty$. We will also assume a moment bound on the size of the jumps:
  \begin{equation}
  \label{ass2}
  \sup_{\bx \in \Sigma} \Exp [ \| \xi_{t+1} - \xi_t \|^\gamma \mid \xi_t = \bx ] < \infty .
  \end{equation}
  
  By an analysis (presented in Section \ref{rwproof})
  of the process $X$ defined by $X_t = \| \xi_t \|$,
  we will see that the following result is a consequence
  of our general Theorems \ref{trans} and  \ref{llnthm}. The condition
  $\limsup_{t \to \infty} \| \xi_t \| = \infty$ a.s.~is
  ensured by, for instance,
   a reasonable `irreducibility' condition,
  such as (A1) in \cite{mmw} in the case where $\Sigma = \Z^d$.
 
 \begin{theo}
 \label{rwthm}
 Suppose that for some $\rho \in (0,\infty)$ and $\beta \in (0,1)$, (\ref{ass1}) holds and that
 $\limsup_{t \to \infty} \| \xi_t \| = \infty$ a.s.. Then 
 \begin{itemize}
 \item [(i)] if (\ref{ass2}) holds for some $\gamma > 1+\beta$, $\| \xi_t \| \to \infty$ a.s.~as $t \to \infty$;
 \item[(ii)] if (\ref{ass2}) holds for some $\gamma > 2+2\beta$,  
 \[ \frac{\| \xi_t \|}{t^{1/(1+\beta)}} \longrightarrow \lambda (\rho, \beta), ~{\rm a.s.}, \]
 as $t \to \infty$,
 for $\lambda$ the constant defined at (\ref{ddef}).
 \end{itemize}
 \end{theo}
 
 The {\em qualitative} behaviour of non-homogeneous random walks $\Xi$ satisfying
 (\ref{ass1}) and (\ref{ass2}) of this section
 was studied
 under slightly different conditions in \cite{mmw}
 (it was assumed there that $\Sigma = \Z^d$ and that the
 moment bound (\ref{ass2}) hold for $\gamma =2$).
 In particular,
 if we impose an additional
 condition on the non-radial drift field such as
 \[ \sup_{\bu \in \R^d : \bu \cdot \bx = 0, \, \| \bu \| =1} \left| \Exp 
 [ (\xi_{t+1} - \xi_t ) \cdot \bu \mid \xi_t = \bx ] \right | = O ( \| \bx \|^{-\beta -\eps} ) ,\]
 for some $\eps >0$,
 it was shown in \cite[Theorem 2.2]{mmw} that $\Xi$ has 
 a limiting direction. That is, there exists a random unit vector $\bu$ such that
 $\xi_t / \| \xi_t \| \to \bu$ a.s., as $t \to \infty$. Combined with
 our quantitative result, Theorem \ref{rwthm}(ii), this implies that
  under the conditions
 of Theorem 2.2 of \cite{mmw}, with the stronger
 condition that (\ref{ass2}) holds for $\gamma > 2+2\beta$,
 \[  \frac{ \xi_t }{t^{1/(1+\beta)}} \longrightarrow \bu \lambda (\rho, \beta), ~{\rm a.s.}, ~{\rm as}~ t \to \infty. \]

\section{Proofs of theorems}
\label{sec:proofs}

\subsection{Overview of the proofs}
\label{outline}

In \cite[Section 3]{mvw}, general techniques were developed
for obtaining almost-sure bounds for stochastic processes
using Lyapunov functions. 
In \cite[Section 4]{mvw}
those techniques were applied to the {\em critical}
regime of the Lamperti problem (i.e., drifts of order $1/x$ at $x$).
The results of \cite[Section 3]{mvw} are a useful starting point for us,
enabling us to prove Theorem \ref{thm0},
but they yield bounds that are considerably less sharp than
those that we ultimately require for Theorem \ref{thm1}. The
sharp bounds in Theorem \ref{thm1}, and the second-order
behaviour in Theorem \ref{cltthm}, require a different approach.
Throughout Sections \ref{proof0},
\ref{proof1} and \ref{cltproof}  we work with the process $X_t^{1+\beta}$,
for which we can establish sharp estimates (see Lemma \ref{lem5} below). In particular, Doob's
decomposition for $X_t^{1+\beta}$ will be the basis for
our proofs of Theorems \ref{thm1} and \ref{cltthm}.

The proof of Theorem \ref{trans} (in Section \ref{proof2}) is somewhat different in flavour,
and uses ideas more closely related to those of Lamperti \cite{lamp1}.
The proof of Theorem \ref{rwthm} (in Section \ref{rwproof}) demonstrates the utility
of our general results in dealing with multi-dimensional processes.

\subsection{Proof of Theorem \ref{thm0}}
\label{proof0}

To prove Theorem \ref{thm0}, we
will need the following result, contained in \cite[Theorem 3.2]{mvw}.
Let $(Y_t)_{t \in \Z^+}$ be an $(\F_t)_{t \in \Z^+}$-adapted
process taking values in an unbounded subset of $[0,\infty)$.

\begin{lm}
\label{lem13} Let  $B \in (0, \infty)$ be such that, for all $t \in \Z^+$,
\[  \Exp [ Y_{t+1}-Y_t \mid \F_t ] \leq B ~{\rm a.s.}. \]
Then
for any $\eps>0$, a.s., for all but finitely many $t \in \Z^+$,
 \begin{align*}
 \sup_{0 \leq s \leq t} Y_s  \leq   t (\log t)^{1+\eps}. \end{align*}
\end{lm}

To prove Theorem \ref{thm0}, we will apply
Lemma \ref{lem13} to the process $Y_t = X_t^{1+\beta}$.
To show that this choice of $Y_t$
satisfies the hypothesis of Lemma \ref{lem13}, we thus need
to show that the expected increment is bounded above. Lemma \ref{lem5} below will take care of this.
 First we need a technical estimate.
For ease of notation we write
$\Delta_t := X_{t+1} - X_t$ throughout the remainder of the paper. 

 \begin{lm}
 \label{lem4}
 Suppose that (A0) holds and that (\ref{moments}) holds for 
  $\gamma \in (0,\infty)$.
 For any $r \in (0,\gamma)$ and $\delta \in (0,1)$ there exists
  $C \in (0,\infty)$
  for which, for all $t \in \Z^+$ and all $x \geq 0$,
 \[   \Exp \left[ | \Delta_t |^r \1 \{ | \Delta_t | \geq x^{1-\delta} \} \mid \F_t \right] \leq C (1+x)^{-(\gamma-r)(1-\delta)} , ~{\rm a.s.}. \]
      \end{lm}
  \proof
  Suppose that (\ref{moments}) holds and fix $\delta \in (0,1)$. 
 By Markov's inequality,
\begin{equation}
\label{eq6}
\Pr [ | \Delta_t | \geq x^{1-\delta} \mid \F_t ] \leq
 x^{\gamma(\delta-1)} \Exp [  | \Delta_t |^\gamma \mid \F_t ]
= O (x^{\gamma (\delta-1)}) , ~{\rm a.s.},
\end{equation}
 using (\ref{moments}).
   Now
  by H\"older's inequality, for any $r \in (0,\gamma)$,  
  \begin{align*} \Exp \left[ | \Delta_t |^{r}  \1 \{ | \Delta_t | \geq
  x^{1-\delta } \} \mid \F_t \right]
  & \leq
  \Exp  [ | \Delta_t |^{\gamma} \mid \F_t ]^{\frac{r}{\gamma}}
  \Pr [ | \Delta_t | \geq
  x^{1-\delta } \mid \F_t ]^{1 - \frac{r}{\gamma}} \\
  & = O ( x^{(\gamma-r) (\delta -1) } ), ~{\rm a.s.},
  \end{align*}
  using (\ref{moments}) and (\ref{eq6}).  
\qed\\

The next result gives (in part (i)) the desired
upper bound for the expected increments of $X_t^{1+\beta}$, and also
provides (in part (ii)) a corresponding lower bound.
Part (iii) is a technical estimate on the higher
moments of the increments that we will need later in our
proof of Theorem \ref{thm1}.

 \begin{lm}
 \label{lem5}
 Suppose that  (A0) holds and that for   $\beta \in [0,1)$,
 $\gamma > 1+\beta$, (\ref{moments}) holds.
 \begin{itemize}
 \item[(i)]
  Suppose that for some $A \in (0,\infty)$, $\limsup_{x \to \infty} ( x^\beta \omu_1 (x) ) \leq A$.
  Then for any $\eps >0$ there exists $K \in (0,\infty)$ such that for all $t \in \Z^+$,
  on $\{ X_t > K \}$,
     \begin{align}
   \label{mom1}
   \Exp [ X_{t+1}^{1+\beta}
  - X_t^{1+\beta}
   \mid \F_t ]  \leq A (1+\beta) + \eps, ~{\rm a.s.}.
   \end{align}
   \item[(ii)] Suppose that for some $a \in (0,\infty)$, $\liminf_{x \to \infty} ( x^\beta \umu_1 (x) ) \geq a$.
   Then for any $\eps >0$ there exists $K \in (0,\infty)$ such that for all $t \in \Z^+$,
  on $\{ X_t > K \}$,
   \begin{align}
   \label{mom1low}
  \Exp [ X_{t+1}^{1+\beta}
  - X_t^{1+\beta}
   \mid \F_t ]   \geq a (1+\beta) - \eps, ~{\rm a.s.}.
   \end{align}
    \item[(iii)]
   Let $r \in [1, \frac{\gamma}{1+\beta} )$.
Then there exists $C \in (0,\infty)$ such that for   all $t \in \Z^+$, 
 \begin{equation}
  \label{mom2}
\Exp [ | X^{1+\beta}_{t+1} - X^{1+\beta}_t |^{r} \mid \F_t ]   \leq C X_t^{\beta r}, ~{\rm a.s.}.
\end{equation}
   \end{itemize}
     \end{lm}
  \proof
  In this proof and all of the proofs that follow, $C$ will denote
a constant whose value may change from line to line.
Recall that $\Delta_t = X_{t+1} - X_t$.
  First we prove parts (i) and (ii) of the lemma.  Let $\delta \in (0,1)$ and define
  the event $E_t := \{ | \Delta_t | < X_t^{1-\delta} \}$; denote
  the complement of $E_t$ by $E^\c_t$. 
 The basic idea is as follows. We will show that  the difference $X_{t+1}^{1+\beta} - X_t^{1+\beta}$
 can be written as
 \[ (1+\beta ) X_t^\beta \left( \Delta_t - \Delta_t \1 ( E^\c_t ) + R_1 (X_t, \Delta_t ) \right) + R_2 (X_t, \Delta_t) ,\]
 where $\Exp [ |\Delta_t \1 ( E^\c_t ) | \mid \F_t ] = o(X_t^{-\beta} )$, $\Exp [ |R_1 (X_t, \Delta_t ) | \mid \F_t ] = o(X_t^{-\beta})$
 and $\Exp [ |R_2 (X_t, \Delta_t)| \mid \F_t ] = o(1)$. Then on taking expectations we see that the dominant
 term is $(1+\beta) X_t^\beta \Exp [ \Delta_t \mid \F_t ]$, which gives the results in (i) and (ii). In the above display the
  term $R_1$ comes from the error term in the Taylor expansion of $(X_{t}+\Delta_t)^{1+\beta} - X_t^{1+\beta}$  on the event $E_t$,
  while the term $R_2$ is $X_{t+1}^{1+\beta} - X_t^{1+\beta}$ on the event $E^\c_t$, which is an event of small probability
  under our assumption of (\ref{moments}).
  
  We now give the details of the argument sketched above.
 Since $X_t + \Delta_t \geq 0$,
  Taylor's theorem with Lagrange form for the remainder implies that
  \begin{align*} X_{t+1}^{1+\beta} - X_t^{1+\beta}   = (X_t + \Delta_t)^{1+\beta} - X_t^{1+\beta} 
    = (1 + \beta) X_t^\beta \Delta_t \left( 1 + \eta \frac{\Delta_t}{X_t} \right)^\beta ,\end{align*}
  where $\eta = \eta (X_t, \Delta_t) \in [0,1]$.
   Since there exists $C \in (0,\infty)$ such that
  $| (1+ y)^\beta - 1 | \leq C | y|$ for all $y \in [-1,1]$, we have that
  \[ \Delta_t \left( 1 + \eta \frac{\Delta_t}{X_t} \right)^\beta \1 ( E_t )
  = \Delta_t \1 ( E_t ) + R_1 (X_t, \Delta_t) ,\]
  where $| R_1 (X_t , \Delta_t ) | \leq C | \Delta_t |^2 X_t ^{-1} \1 ( E_t ) \leq C | \Delta_t |^{1+\beta} X_t^{-1+(1-\delta)(1-\beta)}$.
  Since (\ref{moments}) holds for $\gamma > 1+\beta$, it follows that
  \begin{equation}
  \label{eq50}
   \Exp [ | R_1 (X_t, \Delta_t ) | \mid \F_t ] \leq C X_t^{-\beta-\delta(1-\beta)} = o (X_t^{-\beta}) ,\end{equation}
  as $X_t \to \infty$, since $\beta < 1$.
    Then we have that
  \begin{equation}
  \label{xeq}
   X_{t+1}^{1+\beta} - X_t^{1+\beta} = (1 + \beta ) X_t^\beta \Delta_t \1 ( E_t) + (1+\beta) X_t^\beta R_1 (X_t, \Delta_t)
  + R_2(X_t,\Delta_t) ,\end{equation}
  where
    \begin{equation}
  \label{eq51} R_2(X_t,\Delta_t) = (1 + \beta) X_t^\beta \Delta_t \left( 1 + \eta \frac{\Delta_t}{X_t} \right)^\beta \1 ( E^\c_t) .\end{equation}
  Since on $E_t^\c$ we have
  $X_t^\beta \leq | \Delta_t |^{\beta/(1-\delta)}
  = O ( | \Delta_t |^{\beta + \delta} )$ for $\beta \in [0,1)$ and
  $\delta>0$ small enough,
  it follows that, for small enough $\delta$,
  \begin{equation}
  \label{req}
   | R_2 (X_t, \Delta_t ) | \leq C | \Delta_t |^{1 + \beta +   \delta}  \1 ( E^\c_t)  , ~{\rm a.s.}, \end{equation}
  for some constant $C \in (0,\infty)$ not depending on $t$, $X_t$ or $\Delta_t$.
 Then taking expectations in (\ref{req}) and using
  the $r=1+\beta+\delta$ case of Lemma \ref{lem4} we have that for $t \in\Z^+$,
  \begin{equation} 
  \label{eq52}
  \Exp [ | R_2 (X_t, \Delta_t ) | \mid \F_t ]   \leq C (1+X_t)^{-(\gamma - 1-\beta - \delta)(1-\delta) } = o(1) ,\end{equation}
  as $X_t \to \infty$,
   taking  $\delta \in (0,1)$ small enough so that
  $1+\beta+\delta < \gamma$.
  Also we have that
  \begin{align}
   \Exp [ \Delta_t \1 (E_t) \mid \F_t ]
    = \Exp [ \Delta_t \mid \F_t ] - \Exp [\Delta_t \1 (E_t^\c) \mid \F_t ]  
   =
  \Exp [ \Delta_t \mid \F_t  ] + O ( X_t^{-(\gamma -1)(1-\delta)} ),
  \label{eq0}
  \end{align}
by the $r=1$ case of Lemma \ref{lem4}. Since $\gamma > 1+\beta$ we can choose
$\delta < \frac{\gamma -1-\beta}{\gamma -1}$ so that this last
error term is $o(X_t^{-\beta} )$ as $X_t \to \infty$. 
Then taking expectations in (\ref{xeq}) and using (\ref{eq50}), (\ref{eq52}) and (\ref{eq0}), it follows
that, a.s.,  for all $t \in \Z^+$,
\[  \Exp [ X_{t+1}^{1+\beta} - X_t^{1+\beta} \mid \F_t  ] = (1 + \beta) X_t^\beta  \Exp [ \Delta_t \mid \F_t ] + o(1).\]
Under the conditions of part (i) of the lemma we have that
$\Exp [ \Delta_t \mid \F_t ] \leq (A + o(1)) X_t^{-\beta}$,
as $X_t \to \infty$, a.s.. 
Then (\ref{mom1}) follows. Similarly (\ref{mom1low}) follows under the conditions of part (ii) of the lemma.

It remains to prove part (iii) of the lemma. From (\ref{xeq}) and (\ref{req}), together with the
fact that $| X_t^\beta R_1(X_t, \Delta_t ) | \leq C |\Delta_t |^{1+\beta}$, we have that,
 for $\delta>0$,
\[ |  X_{t+1}^{1+\beta} - X_t^{1+\beta} |^{r} \leq C \left( X_t^\beta | \Delta_t | + 
| \Delta_t |^{1 + \beta + \delta} \right)^{r} , ~{\rm a.s.}.\]
Since $r \geq 1$, Minkowski's inequality implies that
\[ \Exp [  |  X_{t+1}^{1+\beta} - X_t^{1+\beta} |^{r} \mid \F_t ]
\leq C \left( X_t^\beta \Exp [ | \Delta_t |^{r} \mid \F_t]^{1/r} + \Exp [ | \Delta_t |^{(1 + \beta + \delta) r} \mid \F_t ]^{1/r} \right)^{r}, ~{\rm a.s.}.\]
Taking $\delta$ small enough so that $(1 + \beta + \delta )r  \leq \gamma$, which
we can do since $r < \gamma / (1+\beta)$,
we have
from (\ref{moments})  that both of the expectations
on the right-hand side of the last display are uniformly bounded above. Thus (\ref{mom2}) follows.
\qed\\
 
We can now give the proof of Theorem \ref{thm0}. \\

\noindent
{\bf Proof of Theorem \ref{thm0}.}
Lemma \ref{lem5}(i)
shows that under the conditions of Theorem \ref{thm0},
it is legitimate to apply
Lemma \ref{lem13} to the process $Y_t = X_t^{1+\beta}$. This
yields the result. \qed

\subsection{Proof of Theorem \ref{thm1}}
\label{proof1}

Armed with the estimates in Lemma \ref{lem5}, we can now work towards a proof
of Theorem \ref{thm1}. The next result shows that for large $t$,
$X_t^{1+\beta}$ is, to first order,
 well-approximated by the quantity $A_t$ defined
by $A_0 := 0$ and for $t \in \N$ by
\begin{equation}
\label{adef}
 A_t := \sum_{s=0}^{t-1}  \Exp [   X^{1+\beta}_{s+1} - X^{1+\beta}_s   \mid \F_s ] .\end{equation}
 Under the conditions of Theorem \ref{thm1}, $A_t$ will be seen to grow
 linearly with $t$.

 \begin{lm}
 \label{lem:decom}
 Suppose that (A0) holds,   that for $\beta \in [0,1)$,
 $\limsup_{x \to \infty} ( x^\beta \omu_1 (x) ) < \infty$, and that
 (\ref{moments}) holds for some $\gamma > 2 + 2\beta$. Define $A_t$ by (\ref{adef}).
 Then as $t \to \infty$,
 \[ t^{-1} \left| X^{1+\beta}_t - A_t
 \right| \to 0, ~{\rm a.s.}.\]
 \end{lm}
 \proof
 Write $Y_t := X^{1+\beta}_t$.
 By Doob's decomposition (see e.g.~\cite[p.~120]{williams}),
 taking $D_t := \Exp [ Y_{t+1} - Y_t \mid \F_t ]$ and writing
 $A_t := \sum_{s=0}^{t-1} D_s$, we have that $(M_t)_{t \in \Z^+}$,
 defined by $M_0 := Y_0$ and for $t \in \N$ by
 \begin{equation}
 \label{decomp}
  M_t := Y_t - A_t = Y_t - \sum_{s=0}^{t-1} D_s, \end{equation}
 is a martingale adapted to $(\F_t)_{t \in \Z^+}$.
Taking expectations in the identity $M_{t+1}^2 - M_t^2 = (M_{t+1} - M_t)^2 + 2 M_t (M_{t+1} - M_t)$,
we see by the martingale property that
 \begin{equation}
 \label{eq1} \Exp [ M_{t+1}^2 - M_t^2 \mid \F_t ]
 =  \Exp [ ( M_{t+1} - M_t )^2 \mid \F_t ] .\end{equation}
 Moreover, by (\ref{decomp}),
 \begin{align}
 \label{eq2} \Exp [ ( M_{t+1} - M_t )^2 \mid \F_t ] & = \Exp [ ( Y_{t+1} - Y_t - D_t )^2 \mid \F_t ] \nonumber\\
& = \Exp [ ( Y_{t+1} - Y_t )^2 \mid \F_t ] - D_t^2,
 \end{align}
 where we have expanded the term $( Y_{t+1} - Y_t - D_t )^2$ and used the fact that $D_t$ is
 $\F_t$-measurable.
Now by (\ref{eq1}), (\ref{eq2}) and the $r =2$ case of (\ref{mom2}) (which
is valid since $\gamma > 2 (1+\beta)$) we have that for all $t \in \N$,
\[ \Exp [ M^2_{t+1} - M_t^2 ] = \Exp  [ \Exp [ M^2_{t+1} - M_t^2  \mid \F_t ]  ] \leq C \Exp [ X_t^{2\beta} ] .\]
Now since $\beta \in [0,1)$, Jensen's inequality implies that
\[ \Exp [ X_t^{2 \beta} ] \leq \left( \Exp [ X_t^{1+\beta} ] \right)^{\frac{2\beta}{1+\beta}} = O \left(t^{\frac{2\beta}{1+\beta}} \right) ,\]
by  (\ref{mom1}) and the fact that $Y_0$ is uniformly bounded (from the final part of (A0)). 
Thus   $M_t^2$ is a nonnegative submartingale with
  \[
  \Exp [ M_t^2 ] \leq \Exp [ Y_0^2 ] + \sum_{s=0}^{t-1} \Exp [   M_{s+1}^2 - M_s  ^2 ]
  = O \left( t^{\frac{1+3\beta}{1+\beta}} \right) ,\]
  again using the fact that $Y_0$ is uniformly bounded.
  Doob's submartingale inequality (see e.g.~\cite[p.~137]{williams}) then implies that for any $\eps>0$,
  \[ \Pr \left[ \sup_{0 \leq s \leq t} M_s^2 > t^{\frac{1+3\beta}{1+\beta} +   \eps} \right]
  \leq t^{-\frac{1+3\beta}{1+\beta}  -  \eps} \Exp [M_t^2]
  = O(  t^{-\eps} ).\]
  Hence the Borel--Cantelli lemma  implies that for any $\eps>0$, a.s.,
  \[ \sup_{0 \leq s \leq 2^m} | M_s | \leq (2^m)^{\frac{1+3\beta}{2+2\beta} +  \eps} ,\]
    for all but finitely many $m \in \Z^+$.
  Since for any $t\in \N$ we have $2^{m(t)} \leq t < 2^{{m(t)}+1}$ for some $m(t) \in \Z^+$, we have that
  for any $\eps>0$,
  a.s., for all but finitely many $t \in \Z^+$,
  \[  \sup_{0 \leq s \leq t} | M_s | \leq \sup_{0 \leq s \leq 2^{m(t)+1}} | M_s |
  \leq (2^{m(t)+1})^{\frac{1+3\beta}{2+2\beta} +  \eps} \leq C t^{\frac{1+3\beta}{2+2\beta} +  \eps},\]
  for some $C \in (0,\infty)$ not depending on $t$.
 Since $\beta < 1$, we may take $\eps$ small enough
 so that $\frac{1+3\beta}{2+2\beta} +  \eps \leq 1-\eps$. Then
  we have that
  $| A_t - Y_t | = O(t^{1-\eps})$ as $t \to \infty$, a.s..
\qed 

\begin{rmk} The decomposition in Lemma \ref{lem:decom} is 
 central to our proof of Theorem \ref{thm1}. Here the behaviour of the
   supercritical case ($\beta < 1$) is very different to that of
   the critical case when $\beta=1$ (see \cite[Section 4]{mvw})
   where, even in the transient case,
   there is no decomposition available into
   a dominant `drift' part (like $A_t$)
   and a smaller `variation' part (like $M_t$).
   Thus proving (particularly lower) bounds in the critical case needs
   a rather different approach: see \cite{mvw}.\end{rmk}

  Now we can complete the proof of Theorem \ref{thm1}.\\

\noindent
{\bf Proof of Theorem \ref{thm1}.}
   Under the conditions of the theorem,
  we have from Lemma \ref{lem5} that (\ref{mom1}) and (\ref{mom1low}) hold.
  Moreover, we know from Theorem \ref{trans}
  that $X_t \to \infty$ as $t \to \infty$, a.s.,
  so that the $\eps$ terms in (\ref{mom1}) and (\ref{mom1low})
  may be taken to be arbitrarily small for all $t$ large enough.
  Hence with $A_t$ defined by (\ref{adef}) we have that
  for any $\eps >0$,  a.s.,
  \[ a (1+\beta) - \eps \leq t^{-1} A_t \leq A (1+\beta) + \eps, \]
  for all but finitely many $t \in \Z^+$. Now from Lemma \ref{lem:decom}
  we have that
$X_t^{1+\beta} = A_t + o (t)$, a.s., so that
 for any $\eps >0$,  a.s.,
  \[ a (1+\beta) - \eps \leq t^{-1} X^{1+\beta}_t \leq A (1+\beta) + \eps, \]
  for all but finitely many $t \in \Z^+$.
  This proves the theorem. \qed

\subsection{Proof of Theorem \ref{cltthm}}
\label{cltproof}

The basic ingredients of the proof of Theorem
\ref{cltthm} are already in place in the
decomposition used in the proof
of Lemma \ref{lem:decom}, but we need 
to revisit some of our earlier estimates and
obtain sharper bounds under the conditions
of Theorem \ref{cltthm}.
First we have the following refinement
of Lemma \ref{lem5} in this case.

\begin{lm}
\label{lem6}
Suppose that (A0) holds and that for $\beta \in [0,1)$,
$\rho \in (0,\infty)$, (\ref{limits2}) holds.
Suppose that (\ref{moments}) holds for $\gamma > 2+2\beta$. Then
as $t \to \infty$,
\begin{equation}
\label{eq10}
   \Exp [ X_{t+1}^{1+\beta} - X_t^{1+\beta} \mid \F_t ]
= \rho (1+\beta)  + o (t^{-\frac{1-\beta}{2+2\beta}} ), ~{\rm a.s.}. \end{equation}
If in addition (\ref{var}) holds for $\sigma^2 \in (0,\infty)$, then as $t \to \infty$,
\begin{equation}
\label{eq11}
 \Exp [ ( X_{t+1}^{1+\beta} - X_t^{1+\beta})^2 \mid \F_t ] =
\sigma^2 (1+\beta)^2 \lambda (\rho, \beta)^{2 \beta} t^{\frac{2\beta}{1+\beta}}
(1 + o(1)) , ~{\rm a.s.}.\end{equation}
\end{lm}
\proof 
We follow a similar argument to the proof of Lemma \ref{lem5}.
We again use the notation $E_t := \{ | \Delta_t | < X_t ^{1-\delta} \}$
and $E^\c_t$ for the complementary event.
We need to obtain better estimates
for the error terms in (\ref{xeq}) than we did in the proof
of Lemma \ref{lem5}. For this reason the $\delta \in (0,1)$ there will
not be arbitrarily small, so we cannot use (\ref{req}).
Instead, with $R_2(X_t, \Delta_t)$ as defined at (\ref{eq51}),
since on $E^\c_t$ we have $X_t^\beta \leq | \Delta_t |^{\beta/(1-\delta)}$ where
$\beta < 1$,
 we have that, a.s.,
\[ | R_2(X_t, \Delta_t ) | \leq C | \Delta_t |^{1+\frac{\beta}{1-\delta}} \1 ( E^\c_t ) ,  \]
so that from Lemma \ref{lem4}, a.s.,
\[ \Exp [ | R_2(X_t, \Delta_t ) | \mid \F_t ]
 \leq C (1+X_t)^{(1+\frac{\beta}{1-\delta} - \gamma) (1 - \delta) }
= C (1+X_t)^{1-\delta + \beta - (1-\delta) \gamma } .\]
Now take $\delta = \frac{1-\beta}{2} -\eps$
for $\eps > 0$ small enough so that $\delta >0$. It follows that provided
$\gamma > 2$ we can take $\eps>0$ small enough so that
$\Exp [ | R_2(X_t, \Delta_t ) | \mid \F_t ]
= o (X_t^{-\frac{1-\beta}{2}})$. 
Next recall (see just above (\ref{eq50}))
that $| X_t^\beta R_1 (X_t, \Delta_t ) | \leq C X_t^{\beta-1} | \Delta_t |^2$.
Since (\ref{moments}) holds for $\gamma > 2$, we can take expectations to obtain
\[ \Exp [ | X_t^\beta R_1(X_t, \Delta_t ) | \mid \F_t ] = O ( X_t^{\beta -1} ) = o (X_t^{-\frac{1-\beta}{2}}) ,\]
as $X_t \to \infty$, since $\beta < 1$.
Moreover, from (\ref{eq0})
we have that, again taking $\eps$ small enough and using the fact that $\gamma >2$, 
\[ \Exp [ \Delta_t \1 (E_t)  \mid \F_t ]
= \Exp [ \Delta_t \mid \F_t ] + o( X_t^{-\beta-\frac{1-\beta}{2}} )
= \rho X_t^{-\beta}   + o( X_t^{-\beta-\frac{1-\beta}{2}} ) ,~{\rm a.s.},\]
by (\ref{limits2}). 
With these sharper bounds, from 
 (\ref{xeq}) and the present choice of $\delta$ we obtain
\[  \Exp [ X_{t+1}^{1+\beta} - X_t^{1+\beta} \mid \F_t ] = 
(1+\beta ) \rho + o(X_t^{-\frac{1-\beta}{2}} ) , ~{\rm a.s.}.\]
Under the conditions of the lemma,
Theorem \ref{llnthm} applies so that $X_t \sim \lambda(\rho,\beta) t^{\frac{1}{1+\beta}}$.
Thus we obtain (\ref{eq10}).
The argument for (\ref{eq11}) is similar, starting by squaring
(\ref{xeq}) and then taking $\delta>0$ small enough, so we omit the details. \qed  

\begin{lm}
\label{lem7}
Suppose that (A0) holds and that for $\beta \in (0,1)$,
$\rho \in (0,\infty)$, (\ref{limits2}) holds.
Suppose that (\ref{moments}) holds for $\gamma > 2+2\beta$ and that
(\ref{var}) holds for $\sigma^2 \in (0,\infty)$.
Then with $M_t$ as defined at (\ref{decomp}), we have that as $t\to \infty$,
\[ t^{-\frac{1+3\beta}{2+2\beta}} M_t \tod Z \sigma \lambda (\rho, \beta)^{\beta} \sqrt{ \frac{(1+\beta)^{3}}{(1+3\beta)} }
   ,\]
   where $Z$ is a standard normal random variable.
\end{lm}
\proof
We will apply a standard martingale central limit theorem. Set $M_{t,s} = t^{-\frac{1+3\beta}{2+2\beta}} (M_s - M_{s-1})$
for $1 \leq s \leq t$. For fixed $t$, $(M_{t,s})_s$ is a martingale difference
sequence with $\Exp [ M_{t,s} \mid \F_{s-1} ] =0$. Moreover,
\begin{align*}
\sum_{s=1}^t \Exp [ M_{t,s}^2 \mid \F_{s-1} ]
& = t^{-\frac{1+3\beta}{1+ \beta}} \sum_{s=1}^t \Exp [ (M_s - M_{s-1})^2 \mid \F_{s-1} ] \nonumber\\
& = t^{-\frac{1+3\beta}{1+ \beta}}  \sum_{s=1}^t \left( \Exp [ (X^{1+\beta}_s - X^{1+\beta}_{s-1})^2 \mid \F_{s-1} ] + O(1) \right), ~{\rm a.s.},
\end{align*}
by (\ref{eq2}), using the fact that $|D_{s-1}| = O(1)$ by (\ref{eq10}). Now applying
(\ref{eq11}) we obtain
\begin{align}
\label{var1}
 \sum_{s=1}^t \Exp [ M_{t,s}^2 \mid \F_{s-1} ] & = t^{-\frac{1+3\beta}{1+ \beta}}  
\sigma^2 (1+\beta)^2 \lambda (\rho, \beta)^{2 \beta} (1+o(1)) \sum_{s=1}^t s^{\frac{2\beta}{1+\beta}} \nonumber\\
& = \sigma^2 \frac{(1+\beta)^3}{1+3\beta} \lambda (\rho, \beta)^{2 \beta} + o(1), ~{\rm a.s.},\end{align}
where we have used the fact that $\beta >0$ to obtain the $o(1)$ bound in the first equality.
We also need to verify a form of the conditional Lindeberg condition. We claim that, for any $\eps>0$,
\begin{equation}
\label{lind}
\sum_{s=1}^t \Exp [ M_{t,s}^2 \1 \{ |M_{t,s}| > \eps \} \mid \F_{s-1} ]  = o(1), ~{\rm a.s.},
\end{equation}
as $t \to \infty$. To see this, take $p \in (2, \frac{\gamma}{1+\beta})$.
By the elementary inequality $|M_{t,s}|^2 \1 \{ |M_{t,s}| > \eps \} \leq \eps^{2-p} | M_{t,s} |^p$
we have that, for any $\eps>0$,
\[ \Exp [ M_{t,s}^2 \1 \{ |M_{t,s}| > \eps \} \mid \F_{s-1} ]
= O ( \Exp [ |M_{t,s} |^p \mid \F_{s-1} ] ) .\]
Then we have that 
\begin{align*}
\Exp [ |M_{t,s} |^p \mid \F_{s-1} ] & = t^{-\frac{(1+3\beta)p}{2+2\beta}} \Exp [ (M_s - M_{s-1})^p \mid \F_{s-1} ] \\
& = t^{-\frac{(1+3\beta)p}{2+2\beta}} \Exp [ (X^{1+\beta}_s - X^{1+\beta}_{s-1} - D_{s-1})^p \mid \F_{s-1} ] ,
\end{align*}
where by (\ref{eq10}), $|D_{s-1}| = O(1)$, and by the
 $r=p$ case of (\ref{mom2}),
 \[  \Exp [ (X^{1+\beta}_s - X^{1+\beta}_{s-1} )^p \mid \F_{s-1} ] \leq C X_{s-1}^{\beta p} \leq C s^{\frac{\beta p}{1+\beta}}, ~{\rm a.s.}, \]
by Theorem \ref{llnthm}. Hence by Minkowski's inequality, 
\[ \Exp [ |M_{t,s} |^p \mid \F_{s-1} ] \leq C  t^{-\frac{(1+3\beta)p}{2+2\beta}} s^{\frac{\beta p}{1+\beta}} , ~{\rm a.s.}, \]
for some $C \in (0,\infty)$.
 Thus we obtain
\[ \sum_{s=1}^t \Exp [ M_{t,s}^2 \1 \{ |M_{t,s}| > \eps \} \mid \F_{s-1} ] 
\leq C t^{1+\frac{\beta p}{1+\beta}-\frac{(1+3\beta)p}{2+2\beta}} , ~{\rm a.s.}.\]
From this last bound we verify (\ref{lind}) since $p>2$.
Given (\ref{var1}) and (\ref{lind}), we can apply
a standard central limit theorem for
martingale differences (e.g.~\cite[Theorem 35.12, p.~476]{billingsley}) to complete the proof. \qed\\

\noindent
{\bf Proof of Theorem \ref{cltthm}.}
Recall the decomposition at (\ref{decomp}). Under the conditions
of the theorem, Theorem \ref{llnthm} and the proof of Lemma
\ref{lem:decom}
imply that $M_t = o(A_t)$, a.s., so that
\[ X_t = (A_t +M_t)^{\frac{1}{1+\beta}} = A_t^{\frac{1}{1+\beta}} + \frac{1}{1+\beta} M_t A_t^{-\frac{\beta}{1+\beta}} (1+o(1))
  .\]
  Here we have from (\ref{eq10}) that, a.s., $A_t = \rho (1+\beta) t + o (t^{\frac{1+3\beta}{2+2\beta}} )$.
  It follows that, a.s.,
  \[ X_t = \lambda(\rho,\beta) t^{\frac{1}{1+\beta}} + o(t^{1/2}) + \frac{1}{1+\beta} \lambda(\rho,\beta)^{-\beta}
  t^{-\frac{\beta}{1+\beta}} M_t (1+o(1)) .\]
  Rearranging we obtain, a.s.,
  \[ \frac{X_t - \lambda(\rho,\beta) t^{\frac{1}{1+\beta}}}{t^{1/2}}
   = \frac{1}{1+\beta} \lambda(\rho,\beta)^{-\beta} t^{-\frac{1+3\beta}{2+2\beta}} M_t (1+o(1)) + o(1).\]
Now on letting $t \to \infty$ Lemma \ref{lem7} completes the proof. \qed

\subsection{Proof of Theorem \ref{trans}}
\label{proof2}

Our proof of Theorem \ref{trans} under the
minimal moments conditions stated in that theorem
requires some delicate analysis in a similar
vein to the estimates in
Section \ref{proof1}. The key
is the following lemma.

\begin{lm}
\label{superm}
Suppose that (A0) and (A1) hold, and
 that there exists $\beta \in [0,1)$ such that
(\ref{moments}) holds for $\gamma > 1+ \beta$ and
\[ \liminf_{x \to \infty}  ( x^\beta \umu_1 (x)  ) > 0 .\]
Then there exist $\nu >0$ and $M_0 \in (0,\infty)$  such that
for all $t\in\Z^+$, on $\{ X_t > M_0 \}$,
\[ \Exp [ (1+X_{t+1})^{-\nu} - (1+X_t)^{-\nu} \mid \F_t ] \leq 0, ~{\rm a.s.}. \]
\end{lm}
\proof
For ease of notation write $W_t := (1+X_t)^{-\nu}$ and, as before,
 $\Delta_t = X_{t+1} - X_t$. Note that our assumption
 on $\umu_1$ implies that for some $c>0$,
 $\Exp [ \Delta_t \mid \F_t  ] \geq c X_t^{-\beta}$ a.s., for all sufficiently
 large $X_t$ and all $t$.
Let $\delta \in (0,1)$.
First, since $(1+x)^{-\nu}$ is a decreasing function of $x \geq0$, for any $x \geq 0$ we have
\begin{equation}
\label{eq7}
 W_{t+1} - W_t \leq ( W_{t+1} - W_t  ) \1 \{ | \Delta_t | < x^{1-\delta} \}
 + ( W_{t+1} - W_t  ) \1 \{  \Delta_t  \leq - x^{1-\delta} \} .\end{equation}
Moreover,  
\begin{align*} W_{t+1} - W_t = (1+X_t)^{-\nu} \left[ \left( 1+ \frac{\Delta_t}{1+X_t} \right)^{-\nu} - 1 \right] .\end{align*}
Hence by Taylor's theorem with Lagrange remainder,
\begin{align*} & ( W_{t+1} - W_t  ) \1 \{ | \Delta_t | < X_t^{1-\delta} \} \\
 &~~ = (1+X_t)^{-\nu} 
\left[- \nu \frac{\Delta_t}{1+X_t} \left( 1 +  \frac{\eta \Delta_t}{1+X_t} \right)^{-1-\nu} \right] 
 \1 \{ | \Delta_t | < X_t^{1-\delta} \} ,\end{align*}
where $\eta = \eta ( X_t, \Delta_t) \in [0,1]$. 
Here we have that, since $|(1+y)^{-1-\nu} -1 |\leq C | y|$ for some $C \in (0,\infty)$ and any $y \in (-1,1)$,
\begin{align*} \Delta_t \left( 1 +  \frac{\eta \Delta_t}{1+X_t} \right)^{-1-\nu} 
\1 \{ | \Delta_t | < X_t^{1-\delta} \} &
= \Delta_t  \1 \{ | \Delta_t | < X_t^{1-\delta} \} \\ &~~
 + 
O ( | \Delta_t |^2 X_t^{-1} \1 \{ | \Delta_t | < X_t^{1-\delta} \} ),\end{align*}
as $X_t \to \infty$. 
Hence as $X_t \to \infty$, a.s.,
\begin{equation}
\label{eq5}
  ( W_{t+1} - W_t  ) \1 \{ | \Delta_t | < X_t^{1-\delta} \} = - ( \nu + o(1) ) X_t^{-1-\nu} \Delta_t \1 \{ | \Delta_t | < X_t^{1-\delta} \} + S( X_t, \Delta_t),\end{equation}
  where $| S( X_t, \Delta_t ) | \leq C | \Delta_t |^2 X_t^{-2-\nu} \1 \{ | \Delta_t | < X_t^{1-\delta} \} $.
  We have that
  \begin{equation}
  \label{eq53}
   \Exp [ | S( X_t, \Delta_t ) | \mid \F_t ] \leq C X_t^{-2-\nu} X_t^{(1-\delta)(1-\beta)} \Exp [ | \Delta_t |^{1+\beta} \mid \F_t ]
  = O ( X_t^{-1-\beta-\nu-\delta (1-\beta) } ) ,\end{equation}
  since (\ref{moments}) holds for $\gamma > 1+\beta$.
Moreover, since $\gamma > 1+\beta$, (\ref{eq0}) implies that we can take $\delta>0$ small enough so that, as $X_t \to \infty$,
\begin{equation}
\label{eq4}
 \Exp [ \Delta_t \1 \{ | \Delta_t | < X_t^{1-\delta} \} \mid \F_t ] = \Exp [ \Delta_t \mid \F_t ]
+o  (X_t^{-\beta}), ~{\rm a.s.} .\end{equation}
 Hence taking expectations in (\ref{eq5}), and using (\ref{eq53}) and
 (\ref{eq4})  together with the assumption that
$\Exp [ \Delta_t \mid \F_t ] \geq (c+o(1)) X_t^{-\beta}$, we have that as $X_t \to \infty$,
\begin{equation}
\label{eq8}
 \Exp \left[  ( W_{t+1} - W_t  ) \1 \{ | \Delta_t | < X_t^{1-\delta} \} \mid \F_t  \right]
\leq - ( c \nu + o(1) ) X_t^{-1-\beta-\nu} , ~{\rm a.s.}  .\end{equation}
On the other hand, since $W_t \in [0,1]$ a.s., we have that
\begin{equation}
\label{eq9}
 \Exp \left[ ( W_{t+1} - W_t  ) \1 \{   \Delta_t   \leq - X_t^{1-\delta} \} \mid \F_t  \right]
\leq \Pr [ | \Delta_t | \geq X_t^{1-\delta}  \mid \F_t ] = O (X_t^{\gamma (\delta-1)}) ,\end{equation}
by (\ref{eq6}). This last bound is
$O( X_t^{-1-\beta-\delta})$ provided $\delta \leq (\gamma -1 - \beta)/(1+\gamma)$. From (\ref{eq7})
with (\ref{eq8}) and (\ref{eq9}), we therefore conclude that, a.s., as $X_t \to \infty$,
\[ \Exp [ W_{t+1} - W_t \mid \F_t  ] \leq - ( c \nu + o(1) ) X_t^{-1-\beta-\nu} + O(X_t^{-1-\beta-\delta }) .\]
Now taking $\nu \in (0, \delta)$ completes the proof. \qed\\
 
\noindent
{\bf Proof of Theorem \ref{trans}.}
To complete the proof we use a well-known
martingale idea  (see e.g.~\cite[Theorem 2.2.2]{FMM} in the
countable Markov chain case).
With $M_0$ the constant in Lemma \ref{superm}, let $M > M_0$.
 For $s \in \Z^+$ set $T_s := \min \{t > s: X_t \leq M  \}$,
an $(\F_t)_{t\in \Z^+}$-stopping time.
We proceed to show that, for some $c>0$, on $\{ X_s > 2M\}$,
\begin{equation}
\label{claim2}
\Pr [ T_s = \infty \mid \F_s ] > c , ~{\rm a.s.},
\end{equation}
for all $s \in \Z^+$.
By Lemma \ref{superm},
we have that $(1+X_{t \wedge T_s})^{-\nu}$ is a nonnegative supermartingale
adapted to $(\F_t)_{t \in \Z^+}$ for $t \geq s$. Hence for any given $s$,  $(1+X_{t \wedge T_s})^{-\nu}$
converges a.s.~as $t \to \infty$ to some limit, say $L$.
Then if $X_s  > 2 M$   we have
\[   (1 + 2M)^{-\nu} \geq  (1+X_s)^{-\nu} \geq \Exp [ L \mid \F_s ] ,\]
by the supermartingale property. Moreover, on $\{ T_s < \infty \}$ we have
$(1 + X_{t \wedge T_s} )^{-\nu}$ converges to $(1+ X_{T_s} )^{-\nu}$, so that
\[ \Exp [ L \mid \F_s ]
\geq \Exp [ (1+ X_{T_s} )^{-\nu}  \1 \{T_s < \infty \}
 \mid \F_s ]
\geq (1+M)^{-\nu} \Pr [ T_s < \infty \mid \F_s ] ,\]
since $X_{T_s} \leq M$ a.s..
Thus we obtain
\[ \Pr [ T_s < \infty \mid \F_s ] \leq \left( \frac{1 +2M}{1+M} \right)^{-\nu} <1 -c, \] 
for some $c>0$, and so we obtain (\ref{claim2}), as required.
From the assumption that $\limsup_{t \to \infty} X_t = \infty$ a.s.,
we have that a.s.~there exist infinitely many 
$\Z^+$-valued stopping times $s_1 < s_2 < \cdots$ such that
$X_{s_i} > 2M$. By standard arguments
(such as L\'evy's extension of the Borel--Cantelli lemmas)
we can then conclude from (\ref{claim2}) that
$X_t > M$ for all but finitely many $t \in \Z^+$, a.s..
This argument holds for any $M > M_0$, and
so we have that $\lim_{t\to \infty} X_t = \infty$ a.s.,
completing the proof of transience.
 \qed

\subsection{Proof of Theorem \ref{rwthm}}
\label{rwproof}

Let $\Xi$ be as defined in Section \ref{exam2}, and take
$\F_t = \sigma (\xi_0, \xi_1, \ldots, \xi_t)$.
We will consider the process $X$ defined
by $X_t = \| \xi_t \|$. Thus we are in the
final case described in Section
\ref{tech}, where $X_t$ is a function of a Markov process.
We will show that under the conditions
of Theorem \ref{rwthm}, 
$X$ so-defined is an instance
of the supercritical Lamperti problem and
satisfies the conditions
of Theorem \ref{trans} or \ref{llnthm} as appropriate. 
Write $\S = \cup_{\bx \in \Sigma} \{ \| \bx \| \}$
for the state-space of $X$. The next lemma will allow us
to apply our general theorems with $X_t = \| \xi_t \|$.

\begin{lm}
\label{rwlem}
Suppose that for some $\beta \in (0,1)$, $\rho \in (0,\infty)$ and $\gamma >1+\beta$,
 $\Xi$ satisfies (\ref{ass1}) and (\ref{ass2}). Then
$X$ defined by $X_t = \| \xi_t\|$ is a stochastic process
on $\S$ satisfying
\begin{align}
\label{rw1}
 \sup_{t \in \Z^+} \sup_{\bx \in \Sigma} \Exp [ | X_{t+1} - X_t |^\gamma \mid \xi_t = \bx ] < \infty, \\
\label{rw2}
 \lim_{x \to \infty} (x^{\beta} \umu_1 (x) )  =  \lim_{x \to \infty} ( x^{\beta} \omu_1 (x) )  = \rho .\end{align}
\end{lm}
\proof
For ease of notation write $D_t = \xi_{t+1} - \xi_t$. By the triangle inequality,
\begin{equation}
\label{rw0}
 |  X_{t+1} - X_t |  =  | \| \xi_t + D_t \| - \| \xi_t \| | \leq \| D_t \| .\end{equation}
Thus with (\ref{rw0}), (\ref{ass2}) implies (\ref{rw1}). Thus it remains to
prove (\ref{rw2}).
In this case it suffices to show that as $\| \bx \| \to \infty$
\[
\| \bx \|^{\beta} \Exp [ X_{t+1} - X_t \mid \xi_t = \bx ] \to \rho .\]
Suppose $\xi_t = \bx \in \Sigma$,
 and take $\delta \in (0,1)$.
 Then by (\ref{rw0})
 \[ \Exp [  ( X_{t+1} - X_t ) \1 \{ \| D_t \| > \| \bx \|^{1-\delta} \} ]
 \leq \Exp [  \| D_t \|  \1 \{ \| D_t \| > \| \bx \|^{1-\delta} \} ] 
 = O ( \|\bx\|^{-(\gamma -1 )(1-\delta)} ) ,\]
 by an argument similar to the proof of Lemma \ref{lem4}.
 Since $\gamma > 1+\beta$, we can take $\delta >0$ small
 enough so that this last bound is $o (\|\bx\|^{-\beta})$. On the other hand,
 applying Taylor's theorem on $\R^d$ we have that when $\xi_t = \bx$
\begin{align*}
 ( X_{t+1} - X_t ) \1 \{ \| D_t \| \leq \| \bx \|^{1-\delta} \} & = ( \| \bx + D_t \| - \| \bx \| ) \1 \{ \| D_t \| \leq \| \bx \|^{1-\delta} \} \\
 & =
\| \bx + \eta D_t \|^{-1} ( \eta \|D_t\|^2 + D_t \cdot \bx ) \1 \{ \| D_t \| \leq \| \bx \|^{1-\delta} \}
 ,\end{align*}
where $\eta = \eta ( \bx, D_t ) \in [0,1]$. Hence
\begin{align}
\label{ff1} & ~~ \Exp [  ( X_{t+1} - X_t ) \1 \{ \| D_t \| \leq \| \bx \|^{1-\delta} \} \mid \xi_t = \bx ] 
 \nonumber\\
 & = (1+o(1)) \Exp [ ( D_t \cdot \hat \bx ) \1 \{ \| D_t \| \leq \| \bx \|^{1-\delta} \} \mid \xi_t = \bx ]
\nonumber\\ &~~ + O ( \| \bx \|^{-1} ) \Exp [ \| D_t \|^2 \1 \{ \| D_t \| \leq \| \bx \|^{1-\delta} \}  \mid \xi_t = \bx  ] ,\end{align}
as $\| \bx \| \to \infty$. 
Here we have that
\begin{align*} \Exp [ \| D_t \|^2 \1 \{ \| D_t \| \leq \| \bx \|^{1-\delta} \} \mid \xi_t = \bx ] & \leq
\Exp [ \| D_t \|^{1+\beta} \mid \xi_t = \bx ] \| \bx \|^{(1-\beta)(1-\delta)} \\
& = O ( \| \bx \|^{(1-\beta)(1-\delta)} ),\end{align*}
since (\ref{ass2}) holds for $\gamma > 1+\beta$. Thus the second term on the right-hand side
of (\ref{ff1}) is $O( \| \bx \|^{-\beta -\delta (1-\beta)} ) = o( \| \bx \|^{-\beta} )$, since $\beta < 1$ and $\delta >0$.
Moreover, for the first term on the right-hand side of (\ref{ff1}) we have that
\[ \left| \Exp [ ( D_t \cdot \hat \bx ) \1 \{ \| D_t \| > \| \bx \|^{1-\delta} \}  \mid \xi_t = \bx ] \right|
\leq \Exp [ \| D_t \| \1 \{ \| D_t \| > \| \bx \|^{1-\delta} \}  \mid \xi_t = \bx  ] ,\]
which is $o ( \|\bx\|^{-\beta})$ 
as we saw above. Combining our calculations, we have shown that
\[ \Exp [ X_{t+1} - X_t  \mid \xi_t = \bx ]
= (1+o(1)) \Exp [ (\xi_{t+1} - \xi_t ) \cdot \hat \bx \mid \xi_t = \bx ] + o ( \| \bx \|^{-\beta} ) .\]
Hence from (\ref{ass1}) we obtain (\ref{rw2}). \qed\\
 
  \noindent
{\bf Proof of Theorem \ref{rwthm}.}
Lemma \ref{rwlem} shows
that under the conditions of Theorem \ref{rwthm},
$X$ defined by $X_t = \| \xi_t\|$ satisfies
all the conditions of Theorem \ref{trans},
which yields part (i) of the theorem.
Lemma \ref{rwlem} also shows that,
 provided (\ref{ass2}) holds
for $\gamma > 2+2\beta$, all the conditions
of Theorem \ref{llnthm} are satisfied, which
yields part (ii). \qed

\section*{Acknowledgements}

Some of this work was done while AW was at the University of Bristol,
supported by the Heilbronn Institute for Mathematical Research. The authors
are grateful to the anonymous referees for their careful reading of the paper.








\end{document}